\documentclass[numreferences]{kluwer}

\newcommand{\rset}{\mathbb{R}}

\newcommand{\rmax}{\rset_{\max}}
\newcommand{\rmin}{\rset_{\min}}

\newcommand{\pcle}{\preceq}
\newcommand{\Mat}{\mathrm{Mat}}
\newcommand{\0}{\bf{0}}
\newcommand{\1}{\bf{1}}

\begin{document}

\begin{article}

\begin{opening}

\title{A unifying approach to software and hardware design for scientific
calculations and idempotent mathematics
\thanks{The work was supported by the RFBR grant 99--01--01198.\\
Submitted to \textit{Reliable Computing}.}
}

\author{G. L. \surname{Litvinov}\email{litvinov@islc.msk.su}}
\institute{International Sophus Lie Centre}
\author{V. P. \surname{Maslov}\email{maslov@ipmnet.ru}}
\institute{Moscow State Institute of Electronics and Mathematics}
\author{A. Ya. \surname{Rodionov}\email{bs-rt@east.ru}}
\institute{}

\runningtitle{Scientific calculations and idempotent mathematics}

\begin{ao}
INTERNATIONAL SOPHUS LIE CENTRE\\
Nagornaya, 27--4--72, Moscow 113186 Russia\\
e-mail: LITVINOV@ISLC.MSK.SU
\end{ao}

\begin{abstract}
A unifying approach to software and hardware design generated by ideas of
Idempotent Mathematics is discussed. The so-called idempotent correspondence
principle for algorithms, programs and hardware units is described.
A software project based on this approach is presented.
\end{abstract}


\end{opening}

\section{Introduction}

Numerical computations are still very important in computer applications.
But until recently there has been a discrepancy between numerical methods
and software/hardware tools for scientific calculations. In particular,
numerical programming was not much influenced by the progress in
Mathematics, programming languages and technology. Modern tools for
numerical calculations are not unified, standardized and reliable enough.
It is difficult to ensure the necessary accuracy and safety of calculations
without loss of the efficiency and speed of data processing. It is
difficult to get correct and exact estimations of calculation errors. For
example, standard methods of interval arithmetic~\cite{2} do not allow to
take into account error auto-correction effects~\cite{19} and, as a
result, to estimate calculation errors accurately.

However, new ideas in Mathematics and Computer Science lead to a very
promising approach (initially presented in \cite{20}--\cite{22}). An
essential aspect of this approach is developing a system of algorithms,
utilities and programs based on a new mathematical calculus which is called
{\it Idempotent Analysis}, {\it Idempotent Calculus}, or {\it Idempotent
Mathematics}. For many problems in optimization and mathematical
modeling this Idempotent Analysis plays the same unifying role as
Functional Analysis in Mathematical Physics, see, e.g., \cite{14},
\cite{17}, \cite{28}--\cite{30} and surveys \cite{15}, \cite{21}.

Idempotent Analysis is based on replacing the usual arithmetic operations
by a new set of basic operations (such as maximum or minimum). There is a
lot of such new arithmetics, which are associated with sufficiently rich
algebraic structures called idempotent semirings. It is very important that
many problems that are nonlinear in the usual sense become linear with
respect to an appropriate new arithmetic, i.e., linear over a suitable
semiring (the so-called {\it idempotent superposition principle} \cite{26},
\cite{27}, \cite{17}, which is a natural analog of the well-known
superposition principle in Quantum Mechanics). This `linearity'
considerably simplifies explicit construction of their solutions. Examples
are the Bellman equation and its generalizations, the Hamilton--Jacobi
equation, etc. The idempotent analysis serves as a powerful heuristic tool
to construct new algorithms and apply unexpected analogies and ideas
borrowed, e.g., from mathematical physics and quantum mechanics.

The abstract theory is well advanced and includes, in particular, a new
integration theory, linear algebra and spectral theory, idempotent
functional analysis, idempotent Fourier transforms, etc. Its applications
include various optimization problems such as multi-criteria decision
making, optimization on graphs, discrete optimization with a large
parameter (asymptotic problems), optimal design of computer systems and
computer media, optimal organization of parallel data processing, dynamic
programming, applications to differential equations, numerical analysis,
discrete event systems, computer science, discrete mathematics,
mathematical logic, etc. (see, e.g., \cite{1}, \cite{3},
\cite{5}--\cite{15}, \cite{17}, \cite{21}, \cite{26}--\cite{32} and
references therein).

It is possible to implement this new approach to scientific and numeric
calculations in the form of a powerful software system based on unified
algorithms. This approach ensures the arbitrary necessary accuracy and
safety of numerical calculations and working time reduction for
developers and testers of algorithms because of software unification.

Our approach uses techniques of object oriented and functional
programming (see, e.g., \cite{25}, \cite{16}), which are very convenient
for the design of our (suggested) software system.  Computer algebra
techniques \cite{4} are also used. Modern techniques of systolic
processors and VLSI realizations of numerical algorithms including parallel
algorithms of linear algebra (see, e.g., \cite{18}, \cite{31}) are
convenient for effective implementations of the proposed approach to
hardware design.

There is a regular method based on the idempotent theory for constructing
back-end processors and technical devices intended for realization of
basic algorithms of idempotent calculus and mathematics of semirings. These
hardware facilities can increase the speed of data processing.

\section{Mathematical objects and their computer representations}

Numerical algorithms are combinations of basic operations. Usually these
basic operations deal with `numbers'. In fact these `numbers' are thought
of as members of some numerical {\it domains} (real numbers, integers
etc.). But every computer calculation deals with finite {\it models} (or
finite {\it computer representations}) of these numerical domains. For
example, integers can be modeled by integers modulo $2^n$, real numbers can
be represented by rational numbers or floating-point numbers etc.
Discrepancies between mathematical objects (e.g., `ideal' numbers) and their
computer models (representations) lead to calculation errors.

Due to imprecision of sources of input data in real-world problems, the
data usually come in the form of confidence intervals or other number sets
rather than exact quantities. Interval Analysis (see, e.g.,~\cite{2})
extends operations of traditional calculus from numbers to
number intervals to make possible processing such imprecise data and
controlling rounding errors in computational mathematics.

However, there is no universal model that would be good in all cases, so we
have to use varieties of computer models. For example, real numbers can be
represented by the following computer models:
\begin{itemize}
\item standard floating-point numbers,
\item double precision floating-point numbers,
\item arbitrary precision floating-point numbers,
\item rational numbers,
\item finite precision rational numbers,
\item floating-slash and fixed-slash rational numbers,
\item interval numbers,
\item and others.
\end{itemize}

When examining an algorithm it is often useful to have a possibility to
change computer representations of input/output data. For this aim the
corresponding algorithm (and its software implementation) must be universal
enough.

\section{Universal algorithms}

It is very important that many algorithms do not depend on particular
models of a numerical domain and even on the domain itself. Algorithms of
linear algebra (matrix multiplication, Gauss elimination etc.) are good
examples of algorithms of this type.

Of course, one algorithm may be more universal than another algorithm
solving the same problem. For example, numerical integration algorithms
based on the Gauss--Jacobi quadrature formulas actually depend on computer
models because they use finite precision constants. On the contrary, the
rectangular formula and the trapezoid rule do not depend on models and in
principle can be used even in the case of idempotent integration (see
below).

The so-called object oriented software tools and programming languages
(like C$^{++}$ and Java, see, e.g., \cite{25}) are very convenient for
computer implementation of universal algorithms.

In fact there are no reasons to restrict ourselves to numerical domains
only. Actually it may be a ring of polynomials, a field of rational
functions, or an idempotent semiring. The case of idempotent semirings is
extremely important because of numerous applications.

\section{Idempotent correspondence principle}

There is a nontrivial analogy between Mathematics of semirings and Quantum
Mechanics. For example, the field of real numbers can be treated as a
`quantum object' with respect to idempotent semirings, which in turn can be
treated as `classical' or `semi-classical' objects with respect to the
former.

Let $\rset$ be the field of real numbers and $\rset_+$ the subset
of all non-negative numbers. Consider the following change of variables:
$$
   u \mapsto w = h \ln u,
$$
where $u \in \rset_+ \setminus \{0\}$, $h > 0$; thus $u = e^{w/h}$, $w \in
\rset$. Denote by $\0$ the additional element $-\infty$ and by $S$ the
extended real line $\rset \cup \{\0\}$. The above change of variables has a
natural extension $D_h$ to the whole $S$ by $D_h(0) = \0$; also, we
denote $D_h(1) = 0 = \1$.

Denote by $S_h$ the set $S$ equipped with the two operations $\oplus_h$
(generalized addition) and $\odot_h$ (generalized multiplication) such that
$D_h$ is a homomorphism of $\{\rset_+, +, \cdot\}$ to $\{S, \oplus_h,
\odot_h\}$. This means that $D_h(u_1 + u_2) = D_h(u_1) \oplus_h D_h(u_2)$
and $D_h(u_1 \cdot u_2) = D_h(u_1) \odot_h D_h(u_2)$, i.e., $w_1 \odot_h
w_2 = w_1 + w_2$ and $w_1 \oplus_h w_2 = h \ln (e^{w_1/h} + e^{w_2/h})$.
It is easy to prove that $w_1 \oplus_h w_2 \to \max\{w_1, w_2\}$ as $h \to
0$.

Denote by $\rmax$ the set $S = \rset \cup \{\0\}$ equipped with operations
$\oplus = \max$ and $\odot = +$, where ${\0} = -\infty$, ${\1} = 0$ as above.
Algebraic structures in $\rset_+$ and $S_h$ are isomorphic; therefore
$\rmax$ is a result of a deformation of the structure in $\rset_+$.

We stress the obvious analogy with the quantization procedure, where $h$ is
the analog of the Planck constant. In these terms, $\rset_+$ (or $\rset$)
plays the part of a `quantum object' while $\rmax$ acts as a
`classical' or `semi-classical' object that arises as the result of a
{\it dequantization} of this quantum object.

Likewise, denote by $\rmin$ the set $\rset \cup \{\0\}$ equipped with
operations $\oplus = \min$ and $\odot = +$, where ${\0} = +\infty$ and
${\1} = 0$. Clearly, the corresponding dequantization procedure is
generated by the change of variables $u \mapsto w = -h \ln u$.

Consider also the set $\rset \cup \{\0, \1\}$, where ${\0} = -\infty$,
${\1} = +\infty$, together with the operations $\oplus = \max$ and
$\odot=\min$. Obviously, it can be obtained as a result of a `second
dequantization' of $\rset$ or $\rset_+$.

The algebras presented in this section are the most important examples of
idempotent semirings, the central algebraic structure of Idempotent
Analysis.

Consider a set $S$ equipped with two algebraic operations: {\it addition}
$\oplus$ and {\it multiplication} $\odot$. The triple $\{S, \oplus,
\odot\}$ is an {\it idempotent semiring} if it satisfies the following
conditions (here and below, the symbol $\star$ denotes any of the two
operations $\oplus$, $\odot$):
\begin{itemize}
\item the addition $\oplus$ and the multiplication $\odot$ are associative:
$ x \star (y \star z) = (x \star y) \star z$ for all $x, y, z\in S$;
\item the addition $\oplus$ is commutative: $x \oplus y = y \oplus x$ for
all $x,y \in S$;
\item the addition $\oplus$ is {\it idempotent}: $x \oplus x = x$ for all
$x\in S$;
\item the multiplication $\odot$ is {\it distributive} with respect to the
addition $\oplus$: $x\odot(y\oplus z) = (x\odot y)\oplus(x\odot z)$ and
$(x\oplus y)\odot z = (x\odot z)\oplus(y\odot z)$ for all $x, y, z\in S$.
\end{itemize}

A {\it unity} of a semiring $S$ is an element ${\1}\in S$ such that for all
$x \in S$
$$
   {\1} \odot x = x \odot {\1} = x.
$$

A {\it zero} of a semiring $S$ is an element ${\0} \in S$ such that $\0 \neq
\1$ and for all $x \in S$
$$
 {\0}\oplus x = x,\qquad {\0}\odot x = x\odot {\0} = {\0}.
$$

A semiring $S$ is said to be {\it commutative} if $x\odot y=y\odot x$
for all $x,y\in S$. A commutative semiring is called a {\it
semifield} if every nonzero element of this semiring is invertible.
It is clear that $\rset_{\max}$ and $\rset_{\min}$ are semifields.

Note that different versions of this axiomatics are used, see, e.g.,
\cite{1}, \cite{3}, \cite{5}, \cite{6}, \cite{12}, \cite{13}--\cite{15},
\cite{17}, \cite{21}, \cite{23}, \cite{30} and some literature indicated in
these books and papers. Many nontrivial examples of idempotent semirings
can be found, e.g., in \cite{1}, \cite{5}, \cite{6}, \cite{12}, \cite{13},
\cite{14}, \cite{17}, \cite{21}, \cite{23}, \cite{24}, \cite{30}. For
example, every vector lattice or ordered group can be treated as an
idempotent semifield.

The addition $\oplus$ defines the following {\it standard partial order}
on $S$: $x\pcle y$ if and only if $x\oplus y=y$. If $S$ contains a zero
$\0$, then ${\0}\pcle x$ for all $x\in S$. The operations $\oplus$ and
$\odot$ are consistent with this order in the following sense: if
$x\pcle y$, then $x\star z\pcle y\star z$ and $z\star x\pcle z\star y$ for
all $x,y,z \in S$.

The basic object of the traditional calculus is a {\it function}
defined on some set $X$ and taking its values in the field $\rset$
(or $\mathbb{C}$); its idempotent analog is a map $X \to S$, where $X$ is
some set and $S =\rmin$, $\rmax$, or another idempotent semiring. Let us
show that redefinition of basic constructions of traditional calculus in
terms of Idempotent Mathematics can yield interesting and nontrivial results
(see, e.g., \cite{17}, \cite{21}, \cite{23}, \cite{24}, for details,
additional examples and generalizations).

\textsc{Example 1. Integration and measures.} To define an idempotent
analog of the Riemann integral, consider a Riemann sum for a function
$\varphi(x)$, $x \in X = [a,b]$, and substitute semiring operations
$\oplus$ and $\odot$ for operations $+$ and $\cdot$ (usual addition and
multiplication) in its expression (for the sake of being definite, consider
the semiring $\rmax$):
$$
   \sum_{i = 1}^N \varphi(x_i) \cdot \Delta_i \quad\mapsto\quad
   \bigoplus_{i = 1}^N \varphi(x_i) \odot \Delta_i
   = \max_{i = 1, \ldots, N}\, (\varphi(x_i) + \Delta_i),
$$
where $a = x_0 < x_1 < \cdots < x_N = b$, $\Delta_i = x_i - x_{i - 1}$, $i
= 1,\ldots,N$. As $\max_i \Delta_i \to 0$, the integral sum tends to
$$
   \int_X^\oplus \varphi(x)\, dx = \sup_{x \in X} \varphi(x)
$$
for any function $\varphi$:~$X \to \rmax$ that is bounded. In general, for
any set $X$ the set function
$$
   m_\varphi(B) = \sup_{x \in B} \varphi(x), \quad B \subset X,
$$
is called an $\rmax$-{\it measure} on $X$; since $m_\varphi(\bigcup_\alpha
B_\alpha) = \sup_\alpha m_\varphi(B_\alpha)$, this measure is completely
additive. An idempotent integral with respect to this measure is defined as
$$
   \int_X^\oplus \psi(x)\, dm_\varphi
   = \int_X^\oplus \psi(x) \odot \varphi(x)\, dx
   = \sup_{x \in X}\, (\psi(x) + \varphi(x)).
$$

Using the standard partial order it is possible to generalize these
definitions for the case of arbitrary idempotent semirings.

\textsc{Example 2. Fourier--Legendre transform.} Consider the topological
group $G = \rset^n$. The usual Fourier--Laplace transform is defined as
$$
   \varphi(x) \mapsto \widetilde\varphi(\xi)
   = \int_G e^{i\xi \cdot x} \varphi(x)\, dx,
$$
where $\exp(i\xi \cdot x)$ is a {\it character} of the group $G$, i.e.,
a solution of the following functional equation:
$$
   f(x + y) = f(x)f(y).
$$

The idempotent analog of this equation is
$$
   f(x + y) = f(x) \odot f(y) = f(x) + f(y).
$$
Hence `idempotent characters' of the group $G$ are linear functions of the
form $x \mapsto \xi \cdot x = \xi_1 x_1 + \cdots + \xi_n x_n$. Thus the
Fourier--Laplace transform turns into
$$
   \varphi(x) \mapsto \widetilde\varphi(\xi)
   = \int_G^\oplus \xi \cdot x \odot \varphi(x)\, dx
   = \sup_{x \in G}\, (\xi \cdot x + \varphi(x)).
$$
This is the well-known Legendre (or Fenchel) transform.

These examples suggest the following formulation of the idempotent
correspondence principle \cite{20}, \cite{21}:
\begin{quote}
{\it There exists a heuristic correspondence between interesting, useful
and important constructions and results over the field of real (or
complex) numbers and similar constructions and results over idempotent
semirings in the spirit of N. Bohr's correspondence principle in
Quantum Mechanics.}
\end{quote}

So Idempotent Mathematics can be treated as a `classical shadow (or
counterpart)' of the traditional Mathematics over fields.

In particular, an idempotent version of Interval Analysis can be
constructed \cite{24}. The idempotent interval arithmetic appear to be
remarkably simpler than its traditional analog. For example, in the
traditional interval arithmetic multiplication of intervals is not
distributive with respect to interval addition, while idempotent interval
arithmetics conserve distributivity. Idempotent interval arithmetics are
useful for reliable computing.

\section{ Idempotent linearity}

Let $S$ be a commutative idempotent semiring.

The following example of a noncommutative idempotent semiring is very
important.

\textsc{Example 3.} Let $\Mat_n(S)$ be a set of all $S$-valued matrices,
i.e., coefficients of these matrices are elements of $S$. Define the sum
$\oplus$ of matrices $A = \|a_{ij}\|$, $B = \|b_{ij}\| \in \Mat_n(S)$ as $A
\oplus B = \|a_{ij} \oplus b_{ij}\| \in \Mat_n(S)$. The {\it product} of
two matrices $A \in \Mat_n(S)$ and $B \in \Mat_n(S)$ is a matrix $AB \in
\Mat_n(S)$ such that $AB = \|\bigoplus_{k = 1}^m a_{ik} \odot b_{kj}\|$.
The set $\Mat_n(S)$ of square matrices is an idempotent semiring with
respect to these operations. If $\0$ is the zero of $S$, then the matrix $O
= \|o_{ij}\|$, where $o_{ij} = \0$, is the zero of $\Mat_n(S)$; if $\1$ is
the unity of $S$, then the matrix $E = \|\delta_{ij}\|$, where $\delta_{ij}
= \1$ if $i = j$ and $\delta_{ij} = \0$ otherwise, is the unity of
$\Mat_n(S)$.

Now we discuss an idempotent analog of a linear space. A set $V$  is
called a {\it semimodule} over $S$ (or an $S$-semimodule) if it
is equipped with an idempotent commutative associative addition operation
$\oplus_V$ and a multiplication $\odot_V$:~$S \times V \to V$ satisfying the
following conditions: for all $\lambda$, $\mu \in S$, $v$, $w \in V$
\begin{itemize}
\item $(\lambda \odot \mu) \odot_V v = \lambda \odot_V (\mu \odot_V v)$;
\item $\lambda \odot_V (v \oplus_V w)
= (\lambda \odot_V v) \oplus_V (\lambda \odot_V w)$;
\item $(\lambda \oplus \mu) \odot_V v
= (\lambda \odot_V v) \oplus_V (\mu \odot_V v)$.
\end{itemize}
An $S$-semimodule $V$ is called a {\it semimodule with zero} if
${\0} \in S$ and there exists a {\it zero} element
${\0}_V \in V$ such that for all $v \in V$, $\lambda \in S$
\begin{itemize}
\item ${\0}_V \oplus_V v = v$;
\item $\lambda \odot_V {\0}_V = {\0} \odot_V v = {\0}_V$.
\end{itemize}

\textsc{Example 4. Finitely generated free semimodule.} The simplest
$S$-semimodule is the direct product $S^n = \{\, (a_1, \ldots, a_n) \mid
a_j \in S, j = 1, \ldots, n \,\}$. The set of all endomorphisms $S^n \to
S^n$ coincides with the semiring $\Mat_n(S)$ of all $S$-valued matrices
(see Example~3).

The theory of $S$-valued matrices, similar to the well-known
Perron--Fro\-be\-ni\-us theory of nonnegative matrices, is well advanced
and has very many applications, see, e.g., \cite{1}, \cite{3},
\cite{5}--\cite{15}, \cite{17}, \cite{21}, \cite{24}, \cite{29},
\cite{30}--\cite{32}).

\medskip

\textsc{Example 5. Function spaces.} An {\it idempotent function space}
${\cal{F}}(X;S)$ consists of functionals defined on a set $X$ and taking
their values in an idempotent semiring $S$. It is a subset of the set of
all maps $X \to S$ such that if $f(x)$, $g(x) \in {\cal{F}}(X;S)$ and $c
\in S$, then $(f \oplus g)(x) = f(x) \oplus g(x) \in {\cal{F}}(X;S)$ and
$(c \odot f)(x) = c \odot f(x) \in {\cal{F}}(X;S)$; in other words, an
idempotent function space is another example of an $S$-semimodule. If the
semiring $S$ contains a zero element $\0$ and ${\cal{F}}(X;S)$ contains the
zero constant function $o(x) = \0$, then the function space
${\cal{F}}(X;S)$ has the structure of a semimodule with the zero $o(x)$ over
the semiring $S$. If the set $X$ is finite we get the previous example.

Recall that the idempotent addition defines a standard partial order in
$S$. An important example of an idempotent functional space is the space
${\cal{B}}(X;S)$ of all functions $X \to S$ bounded from above with respect
to the partial order $\pcle$ in $S$. There are many interesting spaces of
this type including ${\cal{C}}(X;S)$ (a space of continuous functions
defined on a topological space $X$), analogs of the Sobolev spaces, etc.
(see, e.g., \cite{17}, \cite{21}, \cite{23}, \cite{28}--\cite{30} for
details).

According to the correspondence principle, many important concepts, ideas
and results can be converted from usual Functional Analysis to Idempotent
Analysis. For example, an idempotent scalar product in ${\cal{B}}(X;S)$
can be defined by the formula
$$
   \langle\varphi,\psi\rangle = \int_X^\oplus \varphi(x) \odot \psi(x)\, dx,
$$
where the integral is defined as the `$\sup$' operation (see example 1).

\medskip

\textsc{Example 6. Integral operators.} It is natural to construct
idempotent analogs of integral operators of the form
$$
   K:\,
   \varphi(y) \mapsto (K\varphi)(x)
   = \int_Y^\oplus K(x,y) \odot \varphi(y)\, dy,
$$
where $\varphi(y)$ is an element of a functional space ${\cal{F}}_1(Y;S)$,
$(K\varphi)(x)$ belongs to a space ${\cal{F}}_2(X;S)$ and $K(x,y)$ is a
function $X \times Y \to S$. Such operators are {\it linear}, i.e., they are
homomorphisms of the corresponding functional semimodules. If $S = \rmax$,
then this definition turns into the formula
$$
   (K\varphi)(x) = \sup_{y \in Y}\, (K(x,y) + \varphi(y)).
$$
Formulas of this type are standard in optimization theory.

\section{Superposition principle}

In Quantum Mechanics the superposition principle means that the
Schr\"odi\-n\-ger equation (which is basic for the theory) is linear.
Similarly in Idempotent Mathematics the idempotent superposition principle
means that some important and basic problems and equations (e.g.,
optimization problems, the Bellman equation and its versions and
generalizations, the Hamilton-Jacobi equation) that are nonlinear in the
usual sense can be  treated as linear over appropriate idempotent
semirings, see \cite{26}--\cite{30}, \cite{17}.

Linearity of the Hamilton-Jacobi equation over $\rset_{\min}$ (and
$\rset_{\max}$) can be deduced from the usual linearity (over $\mathbb{C}$)
of the corresponding Schr\"odinger equation by means of the dequantization
procedure described above (in Section 4). In this case the parameter $h$ of
this dequantization coincides with $i\hbar$ , where $\hbar$ is the Planck
constant; so in this case $\hbar$ must take imaginary values (because
$h>0$; see \cite{23} for details). Of course, this is closely related to
variational principles of mechanics.

The situation is similar for the differential Bellman equation, see
\cite{17}.

It is well-known that discrete versions of the Bellman equation can be
treated as linear over appropriate idempotent semirings. The so-called {\it
generalized stationary} (finite dimensional) {\it Bellman equation} has the
form
$$
   X = AX \oplus B,
$$
where $X$, $A$, $B$ are matrices with elements from an idempotent semiring
and the corresponding matrix operations are described in example 3 above;
the matrices $A$ and $B$ are given (specified) and $X$ is unknown.

B.A. Carr\'e \cite{5} used the idempotent linear algebra to show that
different optimization problems for finite graphs can be formulated in
unified manner and reduced to solving these Bellman equations, i.e.,
systems of linear algebraic equations over idempotent semirings. For
example, Bellman's method of solving shortest path problems corresponds to
a version of Jacobi's method for solving systems of linear equations,
whereas Ford's algorithm corresponds to a version of Gauss-Seidel's method.

\section{Correspondence principle for computations}

Of course, the idempotent correspondence principle is valid for
algorithms as well as for their software and hardware implementations
\cite{20}--\cite{22}. Thus:

\begin{quote}
{\it If we have an important and interesting numerical algorithm, then
there is a good chance that its semiring analogs are important and
interesting as well.}
\end{quote}

In particular, according to the superposition principle, analogs of linear
algebra algorithms are especially important. Note that numerical algorithms
for standard infinite-dimensional linear problems over idempotent semirings
(i.e., for problems related to idempotent integration, integral operators
and transformations, the Hamilton-Jacobi and generalized Bellman equations)
deal with the corresponding finite-dimensional (or finite) `linear
approximations'. Nonlinear algorithms often can be approximated by linear
ones. Thus the idempotent linear algebra is a basis for the idempotent
numerical analysis.

Moreover, it is well-known that linear algebra algorithms are convenient
for parallel computations; their idempotent analogs admit parallelization
as well. Thus we obtain a systematic way of applying parallel computation
to optimization problems.

Basic algorithms of linear algebra (such as inner product of two vectors,
matrix addition and multiplication, etc.) often do not depend on concrete
semirings, as well as on the nature of domains containing the elements of
vectors and matrices. Algorithms to construct the closure $A^*={\1}\oplus
A\oplus A^2\oplus\cdots\oplus A^n\oplus\cdots= \bigoplus^{\infty}_{n=1}
A^n$ of an idempotent matrix $A$ can be derived from standard methods for
calculating $({\1} -A)^{-1}$. For the Gauss--Jordan elimination method (via
LU-decomposition) this trick was used in \cite{31}, and the corresponding
algorithm is universal and can be applied both to the Bellman equation and
to computing the inverse of a real (or complex) matrix $({\1} - A)$.
Computation of $A^{-1}$ can be derived from this universal algorithm with
some obvious cosmetic transformations.

Thus it seems reasonable to develop universal algorithms that can deal
equally well with initial data of different domains sharing the same basic
structure \cite{21}, \cite{22}.

\section{Correspondence principle for hardware design}

A systematic application of the correspondence principle to computer
calculations leads to a unifying approach to software and hardware design.

The most important and standard numerical algorithms have many hardware
realizations in the form of technical devices or special processors.  {\it
These devices often can be used as prototypes for new hardware units
generated by substitution of the usual arithmetic operations for its
semiring analogs and by adding a representation of neutral elements $\0$
and} $\1$ (the latter usually is not difficult). Of course, the case of
numerical semirings consisting of real numbers (maybe except neutral
elements) is the most simple and natural \cite{20}--\cite{22}.  Note that
for semifields (including $\rset_{\max}$ and $\rset_{\min}$) the operation
of division is also defined.

Good and efficient technical ideas and decisions can be transposed from
prototypes into new hardware units. Thus the correspondence principle
generates a regular heuristic method for hardware design.  Note that to get
a patent it is necessary to present the so-called `invention formula', that
is to indicate a prototype for the suggested device and the difference
between these devices.

Consider (as a typical example) the most popular and important algorithm
of computing the scalar product of two vectors:
\begin{equation}
   (x, y) = x_1y_1 + x_2y_2 + \dots + x_ny_n.
\end{equation}
The universal version of (1) for any semiring $A$ is obvious:
\begin{equation}
   (x, y) =
   (x_1 \odot y_1) \oplus (x_2 \odot y_2) \oplus \dots \oplus (x_n\odot y_n).
\end{equation}
In the case $A = \rset_{\max}$ this formula turns into the following one:
\begin{equation}
   (x, y) = \max\{ x_1 + y_1, x_2 + y_2, \dots, x_n + y_n\}.
\end{equation}

This calculation is standard for many optimization algorithms, so it is
useful to construct a hardware unit for computing (3). There are many
different devices (and patents) for computing (1) and every such device can
be used as a prototype to construct a new device for computing (3) and even
(2). Many processors for matrix multiplication and for other algorithms of
linear algebra are based on computing scalar products and on the
corresponding `elementary' devices respectively, etc.

There are some methods to make these new devices more universal than their
prototypes. There is a modest collection of possible operations for
standard numerical semirings: max, min, and the usual arithmetic
operations. So, it is easy to construct programmable hardware processors
with variable basic operations. Using modern technologies it is possible to
construct cheap special-purpose multiprocessor chips implementing
reliable, thoroughly tested algorithms. The so-called systolic processors
are especially convenient for this purpose. A systolic array is a
`homogeneous' computing medium consisting of elementary processors, where
the general scheme and processor connections are simple and regular. Every
elementary processor pumps data in and out performing elementary operations
in such a way that the corresponding data flow is kept up in the computing
medium; there is an analogy with the blood circulation, hence
the name `systolic' (see, e.g., \cite{18}, \cite{31}).

Of course, hardware implementations for important and popular basic
algorithms can increase the speed of data processing.

\section{ Correspondence principle for software design}

Software implementations for universal semiring algorithms are not so
efficient as hardware ones (with respect to the computation speed) but are
much more flexible. Program modules can deal with abstract (and variable)
operations and data types. Concrete values for these operations and data
types can be defined by the corresponding input data. In this case concrete
operations and data types are generated by means of additional program
modules. For programs written in this manner it is convenient to use
special techniques of the so-called object oriented (and functional)
design, see, e.g., \cite{25}, \cite{16}. Fortunately, powerful tools
supporting the object-oriented software design have recently appeared
including compilers for real and convenient programming languages (such as
$C^{++}$ and Java).

We propose a project to obtain an implementation of the correspondence
principle approach to scientific calculations in the form of a powerful
software system based on a collection of universal algorithms. This
approach ensures working time reduction for programmers and users because
of software unification. The arbitrary necessary accuracy and safety of
numeric calculations can be ensured as well.

The system contains several levels (including programmer and
user levels) and many modules.

Roughly speaking, it is divided into three parts. The first part
contains modules that implement domain modules (finite representations of
basic mathematical objects). The second part implements universal
(invariant) calculation methods. The third part contains modules
implementing model dependent algorithms. These modules may be
used in user programs written in $C^{++}$ and Java.

The following modules and algorithm implementations are in progress:
\begin{itemize}
\item Domain modules:
\begin{itemize}
   \item infinite precision integers;
   \item rational numbers;
   \item finite precision rational numbers;
   \item finite precision complex rational numbers;
   \item fixed- and floating-slash rational numbers;
   \item complex rational numbers;
   \item arbitrary precision floating-point real numbers;
   \item arbitrary precision complex numbers;
   \item $p$-adic numbers;
   \item interval numbers;
   \item ring of polynomials over different rings;
   \item idempotent semirings $R(\max, \min)$, $R(\max, +)$, $R(\min, +)$,
         interval idempotent semirings
   \item and others.
\end{itemize}
\item Algorithms:
\begin{itemize}
   \item linear algebra;
   \item numerical integration;
   \item roots of polynomials;
   \item spline interpolations and approximations;
   \item rational and polynomial interpolations and approximations;
   \item special functions calculation;
   \item differential equations;
   \item optimization and optimal control;
   \item and others.
\end{itemize}
\end{itemize}

This software system may be especially useful for designers of algorithms,
software engineers, students and mathematicians.

\end{article}

\end{document}